\documentclass[12pt]{article}
\newtheorem{Lemma}{Proposition}

\newcommand{\eps}{\varepsilon}

\newcommand{\dst}{\displaystyle}

\font\msam=msam10 scaled 1200
\font\Bbb=msbm10 scaled 1200
\newcommand{\Rr}{\mbox{\Bbb R}}

\newcommand{\varkappa}{\mbox{\Bbb\symbol{"7B}}}

\newcommand{\WhiteBox}{\phantom{a}\hfill$\mbox{\msam\symbol{3}}$}

\def\rem#1{} 
\def\al{\alpha}
\def\be{\beta}

\def\del{\delta}
\def\ga{\gamma}

\newcommand{\x}{{\bf x}}
\newcommand{\y}{{\bf y}}
\newcommand{\z}{{\bf z}}
\newcommand{\R}{\mbox{\bf R}}

\newcommand{\CH}{{\cal H}}
\newcommand{\CF}{{\cal F}}

\def\tr{\mbox{\rm tr}}
\def\spR{{\rm sp}(2n,\Rr)}
\newcommand{\rk}{{\rm rk}\,}
\begin{document}

\centerline{\bf TWO-SYSTEM OF A HAMILTONIAN SYSTEM}

\vskip 0.5cm
\centerline{\bf M.F.~Kondratieva and S.Yu.~Sadov}

\vskip 1cm
\hangindent=0.5in
\noindent \hspace*{0.5in}
{\small
{\bf Abstract.}
For a Hamiltonian system $\,\dot\x=JH'(\x)\,$ in $\Rr^{2n}\,$
its two-system is defined in the phase space
$\,\Rr^{2n}_{\x}\times \spR_{\Phi}\,$ as follows:
$\;\dot\x=J\,(H'(\x)+\tr(\Phi H''))'$, $\quad
\dot\Phi=[\Phi, H''(\x)J]$.
In a sense, it is a combination of the original system and its system
in variations with feedback. 
We study the Hamiltonian forms of the two-system and
its analogs. In particular, we show a relation
between the Poisson-Lie bracket on the dual space to Lie algebra $\spR$
and the canonical bracket.  
}

\vskip 0.7cm
\subsection*{Introduction}

In this paper we describe one special construction in the theory
of Hamiltonian system called the two-system.
Our point of view is different from that in
\cite{BBK}, where such systems were defined, and in
\cite{BK}, where the term ``two-system" was introduced 
(in the context of semiclassical approximation in quantum mechanics).
Here we treat the two-system as a certain union
of the given Hamiltonian system and its system in variations.
Versions of the construction, their Hamiltonian structures,
and general properties are discussed.

As an example, the following fifth order system of 
ordinary differential equations with parameter $\eps$
$$
\begin{array}{l}
\dst
\dot q =p,\qquad \dot p= -(q +\eps q^3)-3\eps q \al,\\[2ex]
\dst\dot \al =2\be,\qquad \dot \be =-(1+3\eps q^2)\al +\ga,
\qquad
\dot \ga =-2\be(1+3\eps q^2).
\end{array}
$$
is the two-system corresponding to the underlying system with
one degree of freedom and Hamilton's function
$$
 H(q,p)\,=\,\frac{q^2+p^2}{2}\,+\,\eps\,\frac{ q^4}{4}.
$$
This system with $\eps\ll 1$ was analyzed in \cite{Sadov}
by the normal form method. 

\subsection{A hamiltonian view on the system in variations}

\noindent
{\bf a.} $\,$
Consider a Hamiltonian system with $n$ degrees of freedom
\begin{equation}
\label{eq1}
\dot \x =JH'(\x).
\end{equation}
Here $\x=(x_1,\dots, x_{2n})^t$ denotes a point in the phase space
$\Rr^{2n}$ as well as a vector-function $\x(t)$; $H(\x)$ is Hamilton's
function, $H'(\x)$ its gradient
$$
H'(\x)=\left(\frac {\partial H}{\partial x_1},\dots ,\frac {\partial H}
{\partial x_{2n}}\right)^t,
$$
$J$ is the standard symplectic matrix: if
$I_n$ denotes the unity matrix of order $n$, then
$$
  J=\left(\begin{array}{cc} 0& I_n\\-I_n& 0\end{array}\right).
$$

\bigskip
\noindent
{\bf b.} $\,$
Along with system (\ref{eq1}) consider its system in variations
\begin{equation}
\label{eq2}
\dot \y =JH''(\x(t))\y, \qquad
y=(y_1,\dots,y_{2n})^t,
\end{equation}
where $\x(t)$ is a solution of system (\ref{eq1}),
and $H''(\x)$ is the Hesse matrix,
$$
  H''_{ij}=\frac{\partial^2 H}{\partial x_i \partial x_j}.
$$
The following is a simple and well-known fact.

\begin{Lemma}
$\;$
Let $\x(t)$ be any function
{\rm (not necessarily
a solution of system (\ref{eq1}))}.
Then system {\rm(\ref{eq2})}
is a time-dependent Hamiltonian system with respect to variables $\y$,
with quadratic Hamilton's function
\begin{equation}
\label{eq3}
F(\x(t),\,\y) =\frac 1 2 \y^t H''(\x(t))\y.
\end{equation}
\end{Lemma}

\smallskip
\noindent
{\bf c.} $\,$
Our intent is to combine
systems (\ref{eq1}) and (\ref{eq2})
into an autonomous Hamiltonian system.
A straightforward union of (\ref{eq1}) and (\ref{eq2}) suits to this
purpose only in exceptional cases.

\begin{Lemma}
$\;$
Suppose the function $\CF(\x,\y)$ is
such that system {\rm(\ref{eq1})} has the form
\begin{equation}
\label{eq4}
  \dot  \x =J\frac {\partial \CF}{\partial \x}
\end{equation}
and at the same time system {\rm(\ref{eq2})} has the form
\begin{equation}
\label{eq5}
  \dot  \y =J\frac {\partial \CF}{\partial \y}
\end{equation}
Then $H(\x)$ is a polynomial of degree $2$ or less.
\end{Lemma}

\noindent
{\it Proof}. $\,$
It follows from (\ref{eq1}) and (\ref{eq4}) that
$$
  H'(x)=\frac {\partial \CF}{\partial \x}.
$$
Hence $\,\CF(\x,\y)=H(\x) +C(\y)\,$ with some function
$C(\y)$ independent of $\x$. Then (\ref{eq5}) becomes $\dot \y =JC'(\y)$.
It is a system independent of $\x$, which doesn't agree with (\ref{eq2})
unless $H''$ is a constant matrix.
\WhiteBox

\bigskip
Proposition 2 implies that in general any Hamiltonian combination of systems
(\ref{eq1}) and (\ref{eq2})  has to be {\it artificial} in the sense that
it will {\it alter the dynamics}.
Nevertheless such a combination may arise naturally from another
perspective. The Hamiltonian combination we deal with in the
paper, or rather its generalization called the two-system,
first appeared in the context of semiclassical approximation in
quantum mechanics \cite{BBK}. We do not discuss 
quantum mechanical applications here.

\subsection{Hamiltonian combination in the vector form}
\label{vecform}

\noindent
{\it Construction.}
Consider the phase space $\Rr^{4n}$ with coordinates $(\x,\y)$ and a
symplectic structure given by the matrix $J\oplus J$, so that
$(x_i,x_{i+n})$ and  $(y_i,y_{i+n})$ are pairs of canonically conjugate
variables.
Define Hamilton's function in $\Rr^{4n}$ as the sum of
Hamilton's functions of systems (\ref{eq1}) and (\ref{eq2})
\begin{equation}
\label{eq6}
\CH(\x,\y)= H(\x) +\,\frac 1 2\, \y^t H''(\x) \,\y.
\end{equation}
Setting in (\ref{eq4}), (\ref{eq5}) $\,\CF=\CH$, we obtain the Hamiltonian
system
\begin{equation}
\label{eq7}
\begin{array}{rcl}
\dot \x &=& \dst J\,\left( \mathstrut H(\x)+\,\frac 1 2\, \y^t H''(\x)\,\y
\,\right)'\, ,
\\[2ex]
\dot \y &=& \dst J\,H''(\x)\,\y .
\end{array}
\end{equation}
The prime ($\,'\,$) denotes differentiation with respect to $\x$ throughout.

\medskip
\noindent
{\bf Definition 1.}$\,$
We call system (\ref{eq7}) the 
{\it vector form} (for a Hamiltonian combination
of system (\ref{eq1}) and its system in variations (\ref{eq2})$\,$)).

\subsection{Hamiltonian combination in the matrix form}

\noindent
{\bf a.}$\,$
Introduce the $2n\times 2n$ matrix $\,M=\y\y^t\,$ and rewrite
(\ref{eq6}) in the form
 \begin{equation}
 \label{eq7b}
 \CH(\x,\,y)=H(\x)+ \,\frac 1 2 \,\tr ( H''(\x)\,M).
 \end{equation}
The second equation in (\ref{eq7}) is equivalent to
\begin{equation}
\label{eq7a}
\dot M =J\,H''(\x)\, M \,+\,M\,H''(\x)\, J^t.
\end{equation}
Set
$$
\Phi= JM,
\qquad
A=H''J.
$$
Using the identities $\,J^t=-J$, $\,J^2=-I_{2n}$,
write (\ref{eq7a}) in the Lax pair form
$$
\dot \Phi \,=\, [A,\, \Phi]\,\equiv\,A\Phi-\Phi A.
$$
The entire system (\ref{eq7}) becomes
\begin{equation}
\label{eq8}
\begin{array}{rcl}
\dot \x &=& \dst J\,\left(\, H(\x)\,-\,\frac{1}{2}\, \tr(A(\x)\,\Phi)
\,\right)',
\\[2ex]
\dot \Phi &=& \dst  [A(\x),\, \Phi].
\end{array}
\end{equation}

\smallskip
\noindent
{\bf Definition 2.}$\,$
System (\ref{eq8}) with a vector-function $\x(t)$ and a matrix-valued
function $\Phi(t)$ is called the 
{\it matrix form} (for a Hamiltonian combination
of the systems (\ref{eq1}) and (\ref{eq2})$\,$)).

\medskip
\noindent
{\bf Remark}.
The vector form (\ref{eq7}) is contained in the matrix form (\ref{eq8})
and corresponds to the invariant set of such $(\x,\Phi)$
that $M=-J\Phi$ is a symmetric non-negative matrix of rank 1. 
Section \ref{matrank}{\bf c} elaborates on this correspondence.

\medskip
\noindent
{\bf b.}$\,$
An arbitrary matrix $M$ can be decomposed into its symmetric and
antisymmetric parts
\begin{equation}
\label{sadecomp}
M_s=\frac{M+M^t}{2},
\qquad
M_a=\frac{M-M^t}{2}.
\end{equation}
Now take $\,M=-J\Phi$.
Let us show that system (\ref{eq8}) respects the decomposition
(\ref{sadecomp}), and that
only the symmetric part $M_s$ has an effect on the evolution of $\x$.

\begin{Lemma}
$\,$
Let $\,(\x(t),\,JM(t)\,)\,$ be a solution of system
{\rm(\ref{eq8})}, $\;M=M_s+M_a\;$ be the decomposition
{\rm(\ref{sadecomp})}. Denote
$\,\Phi^s=JM_s$, $\,\Phi^a=JM_a$.
Then

\noindent
{\rm(i)} $\,$ $\,(\x(t),\,\Phi^s(t)\,)\,$
is a solution of system {\rm(\ref{eq8})} with $\,\Phi\,$ replaced by
$\,\Phi^s$.

\noindent
{\rm(ii)} $\,$ Evolution of $\,\Phi^a\,$ is given by the equation
$\;\dot\Phi^a=[A,\Phi^a]$.

\noindent
{\rm(iii)} $\,$ If $\,n=1$, then  $\,\dot \Phi^a=0$.
\end{Lemma}

\noindent
{\it Proof}. $\,$
The matrix $M^t$ satisfies the
same equation (\ref{eq7a}) as does $M$
$$
\dot M^t =J H''(x) M^t + M^t H''(x) J^t.
$$
Therefore (\ref{eq7a}) is satisfied separately by
$M_s$ and $M_a$. Correspondingly, the second equation in (\ref{eq8})
is satisfied separately by $\Phi^s$ and $\Phi^a$.

The following calculation shows that the RHS 
of the first equation
in (\ref{eq8}) depends only on $\Phi^s\,$:
$$
\tr \,(A\Phi) =\tr\,(H'' M)\,
\stackrel{*}{=}\,
\tr\,(H'' M_s)=\tr\,(A \Phi^s).
$$
Step $(*)$ holds because
$H''$ is symmetric and
$\;\tr (H''M)=\tr (M^t H'')=\tr(H''M^t)$.
Finally, in case $\,n=1\,$ the matrix $\,\Phi^a\,$ is proportional to
$\,I_2\,$, so $\,\dot \Phi^a=[A,I_2]=0$.
\WhiteBox

\medskip
We are interested in an influence of the described extension
on the original dynamics due to the entanglement of the
$\x$ and $\Phi$ variables. From this point of view, it suffices
to consider system (\ref{eq8}) assuming that matrix $\Phi$
has the form $\Phi=JM$ with a symmetric $M$.
Consequently,  it obeys $\;\Phi^t J +J \Phi =0$.
Recall that $2n\times 2n$ real matrices with this property
are called {\it symplectic}. They form the classical Lie algebra $\,\spR $.

\medskip
\noindent
{\bf Definition 3.} $\,$
System (\ref{eq8}) with $\,\Phi\in\,\spR$
is called the {\it two-system} of Hamilton's system (\ref{eq1}).

\medskip
\noindent
{\bf Remark.} $\,$ The term {\it two-system}
and a family of $n$-systems
were introduced in \cite{BK}. According to \cite{BK},
$n$-system refines the classical limit approximation
for a localized solution of the Schr\"odinger equation
by taking into account all moments of the wave packet up to
order $n$. Higher $n$-systems seem not to be practical because
they add too many degrees of freedom to the original system.
A Hamiltonian formulation
has not been found for $n$-systems with $n>2$.
The situation is quite different for the $2$-system, as
we discuss in the next section.

\subsection{Hamiltonian structure of the two-system}
\label{matrank}

\noindent
{\bf a.} $\,$
By definition, the vector form (\ref{eq7}) is a Hamiltonian
system. However, it is not immediately obvious which
Hamiltonian structure (if any) is related to
to the two-system.
Such a structure was found in \cite{Sadov}.
The result will be recalled in {\bf d}; it
uses the notion of a linear Poisson bracket on $\,\spR$
(more precisely, on its dual space, but we don't make a difference
using the standard inner product on the space of matrices: 
$\,(A,B)=\tr (AB^t)\,$ ).

Preliminary, in {\bf b}, {\bf c} we describe another
Hamiltonian approach to the two-system, considering a family of {\it
multivector forms}, which generalize the vector form.
The phase space of the
two-system is decomposed into invariant subsets corresponding to fixed
ranks and signatures of the symmetric matrix $M=-J\Phi$.

\medskip
\noindent
{\bf b.} $\,$
Let $m_+$, $m_-$ and $m_0$ be the numbers of positive, negative,
and zero eigenvalues of a symmetric $2n\times 2n$ matrix $M$.
Then
$$
 m_+ + m_- + m_0 =2n, \qquad m_+ + m_- = {\rm rk}\,M.
$$
The pair $(m_+,m_-)$ is called {\it signature} of $M$.

\begin{Lemma}
$\,$
Let a symmetric matrix-function $M(t)$ satisfy the equation
$$
J\dot M = [JM, \,A(t)]
$$
with a function $\,A:\,\R\to\spR$. Then
signature of $M(t)$ does not change.
\end{Lemma}

\noindent
{\it Proof.} $\,$
There exists a non-degenerate matrix function $B(t)$ such that
$$
 JM(t)=B(t)\,JM(0)\,B(t)^{-1}
$$
Multiplication by a non-degenerate matrix preserves rank, so
$\,\rk M(t)\,=\,{\rm const}$.
The matrix function $JM(t)$ is continuous and it has
constant rank and real eigenvalues. Nonzero eigenvalues
can not turn to zero, hence their signs don't change.
\WhiteBox

\medskip
\noindent
{\bf Definition 4}. $\,$
Define a subset $\,X(m_+,m_-)\,$ of the phase space $\Rr^{2n}\times
\spR\,$ as the set of all pairs $(\x,\,JM)$ with matrix
$M$ of signature  $(m_+,m_-)$. It is invariant (Proposition 4).
Restriction of the two-system on $\,X(m_+,m_-)\,$ is
called the {\it two-system in signature $\,(m_+,m_-)\,$}.

\medskip
\noindent
{\bf c.}$\,$
{\it Construction} $\,$
(cf.~Sect.~\ref{vecform}).
Let $m_+$, $m_-$ be nonnegative integers and $r=m_+ + m_-\le 2n$.
Take the sum of $r+1$ copies of $\Rr^{2n}$ to be the phase space
of a new dynamical system. Let $\x$, $\;\{\y_i\}=\y_1, \,\dots,\,\y_{m_+}$,
$\;\{\z_j\}=\z_1, \,\dots,\,\z_{m_-}$
be the respective groups of dynamical variables.
Equip the phase space with the standard
symplectic structure $J\oplus J\oplus\cdots \oplus J$ and define
a new Hamilton's function
\begin{equation}
\label{multivec}
\CH(\x,\,\{\y_i\},\,\{\z_j\})\;=\;
  H(\x)\;+\;\sum_{i=1}^{m_+} F(\x,\,\y_i)\; -\;\sum_{j=1}^{m_-}
  F(\x,\,\z_j),
\end{equation}
where $H(\x)$ is Hamilton's function of system (\ref{eq1})
and $F$ is defined in (\ref{eq3}).

\medskip
\noindent
{\bf Definition 5}. $\,$ The system described in the Construction
is called the {\it multi-vector form of signature $(m_+,m_-)$}
(for a Hamiltonian combination of the systems (\ref{eq1}) and
(\ref{eq2})$\,$)).
For brevity, we'll write {\it multivector $(m_+,m_-)$-form}.

\bigskip
\noindent
A symmetric matrix $M$ of signature $(m_+,m_-)$ can be written in the
form
\begin{equation}
\label{mvdecomp}
  M= \sum_{j=1}^{m+}\, \y_j\, \y_j^t
  \,-\,\sum_{i=1}^{m_-} \,\z_i\, \z_i^t
\end{equation}
The matrix built from the vector-columns $\y_i/\|\y_i\|$,
$\z_j/\|\z_j\|$ is determined up to a
pseudo-orthogonal rotation, namely
to an element of the group $O(m_+,m_-)\times I_{m_0}$.  Consequently, the
two-system in signature $\,(m_+,m_-)\,$ can be obtained from the
multivector $(m_+,m_-)$-form by projection. A precise formulation
follows.

\begin{Lemma}
$\,$
Let $\,(\x_0, JM_0)\,$ be the initial point of
a trajectory $\,(\x(t), JM(t)\,)\,$ of the two-system
in the signature $(m_+,m_-)$. Let
$\;\{\y_i(0)\}_1^{m_+}$, $\;\{\z_j(0)\}_1^{m_-}\;$
be the components of some decomposition of $M_0$ of the form
{\rm(\ref{mvdecomp})}.  If $\,(\tilde\x(t)$,
$\,\{\y_i(t)\}$, $\,\{\z_j(t)\}\,)\,$
is the solution of the multivector $\,(m_+,m_-)$-form with initial data
$\,(\x_0$, $\,\{\y_i(0)\}$, $\,\{\z_j(0)\}\,)\,$, then for all $t$
$$
\x(t)=\tilde\x(t)
$$
 and
$$
  M(t)\;=\; \sum_{i=1}^{m_+} \y_i(t) \y_i^t(t)\;-\;
  \sum_{j=1}^{m_-} \z_j(t) \z_j^t(t).
$$
\end{Lemma}

\noindent
{\it Proof.}$\,$
Assuming $M(t)$ in the said form,
it obeys (\ref{eq7a}), if $\y_j$ and $\z_i$ obey (\ref{eq2}).
Also the evolution of $\,\tilde\x(t)\,$ determined by the Hamiltonian
(\ref{multivec}) coincides with that in the first equation (\ref{eq8}).
\WhiteBox

\medskip
Note that a multivector form cannot be subdivided into
vector forms (\ref{eq7})
since $\x(t)$ is subject to a collective effect
by all $\y$'s and $\z$'s.

\medskip
\noindent
{\bf d.}$\,$
The two-system also possesses a non-trivial Hamiltonian structure
with a degenerate Poisson bracket. (Note without going
into details that this fact reflects some general
construction of symplectic geometry \cite[\S~1.5]{KM}.)

Identify the phase space of the two-system with the space
$R^{2n}\times\spR$. For functions on this space,
introduce a Poisson bracket  by the formula
\begin{equation}
\label{eq10}
\{U,V\}(\x,\Phi)= U'\,J\,V'\, +\,2\, \tr
\left(\Phi^t \,[\nabla_{\Phi}U,\nabla_{\Phi}V]\right).
\end{equation}
Here $\,U'\,$ denotes the gradient of $U(\x,\Phi)$
with respect to $x$-variables, and
$\nabla_{\Phi} U$ denotes the matrix from $\spR$
with entries
$$
(\nabla_{\Phi} U)_{ij}=\frac {\partial U}{\partial \Phi_{ij}}.
$$
In particular,
\begin{equation}
\label {eq10b}
\nabla_{\Phi} (\,\det \Phi \,)=(\Phi^{-1})^t\det\Phi;\qquad
\nabla_{\Phi} (\,\tr(\Phi A)\,)=A^t,\quad \forall A.
\end{equation}

\noindent
The bilinear operator (\ref{eq10}) is indeed a
Poisson bracket, i.e.\ it is anti-symmetric and satisfies the Leibnitz and
Jacobi identities.

\begin{Lemma}
The two-system {\rm(\ref{eq8})} can be written in the Hamiltonian form
\begin{equation}
\label{eq11}
x_i =\{x_i,\CH(\x,\Phi)\}, \quad
\dot \Phi_{ij} =\{\Phi_{ij},\CH(\x,\Phi)\}
\end{equation}
with Poisson bracket {\rm(\ref{eq10})} and Hamilton's function
\begin{equation}
\label{eq10a}
 \CH(\x;\Phi)=H(\x)- \frac 1 2 \tr ( A(\x)\Phi), \qquad A(\x)=H''(\x) J
\end{equation}
\end{Lemma}

\noindent
{\it Proof}.
Equivalence of the first equations in (\ref{eq11}) and
(\ref{eq8}) immediately follows from (\ref{eq10}).
Then, $\,\nabla_{\Phi}\Phi_{ij}=E_{ij}\,$, the matrix with the only
nonzero (unity) entry $[ij]$. The second equation in (\ref{eq10b}) gives
$\,\nabla_{\Phi}\CH=-(1/2)\,A^t$. Now (\ref{eq10}) yields
$$
\{\Phi_{ij},\CH(\x,\Phi)\}=\tr\,\left(\Phi^t\,
(-E_{ij}A^t+A^t E_{ij})\right)
=\tr(E_{ij}\,[A,\Phi]^t)=[A,\Phi]_{ij},
$$
so the right-hand sides of the second equations in
(\ref{eq11}) and (\ref{eq8}) agree.
\WhiteBox

\medskip
\noindent
{\bf Remark.}
The Poisson bracket (\ref{eq10}) can be described via the matrix
function $\,\Omega(\x,\Phi)\,$ on the
$\,(2n^2+3n)$-dimensional linear space $\Rr^{2n}\times\spR$.
The matrix $\Omega$ at a generic point has corank
$n$. Moreover, we have the following

\begin{Lemma}
$\,$
Besides Hamilton's function $\CH$, the
two-system has $n$ additional,  
functionally independent integrals, which are Casimir's functions
of the bracket {\rm(\ref{eq10})}. As such integrals, one can take
the coefficiens of the characteristic polynomial $P(\lambda)=\det
(\Phi-\lambda I)$ at $\lambda^0,\,\lambda^2,\dots,\lambda^{2(n-1)}$.
\end{Lemma}

\noindent
{\it Proof}.
Note
that  $\,\det\Phi\,$ is in involution with any function $U(\x,\Phi)$
with respect to the bracket  (\ref{eq10})
$$
  \{\det\Phi,\; U(\x,\Phi)\}=0,
$$
which can be shown using the first of equations (\ref{eq10b}).

Since $\Phi$ is symplectic, its characteristic polynomial is an even
function of degree $2n$. So it has $n$ non-trivial coefficients
corresponding to the powers $\,0,\,2,\dots,2(n-1)\,$ of $\,\lambda$.
They can independently take any values.
\WhiteBox

\smallskip
It is known from the general theory of Hamiltonian systems with degenerate
Poisson brackets (see e.g.\ \cite{Gori}, \cite{DS 1989}) that
the phase space of the 2-system is a disjoint union of invariant
subspaces --- symplectic leafs. Restriction of the bracket on
any symplectic leaf is non-degenerate, hence, at least locally,
there exist canonical coordinates on the leaf, in which the
system has the standard Hamiltonian form.
A point of the phase space in a generic posisiton belongs to
a symplectic leaf of the lowest possible codimension $n$ and the corresponding
nondegenerate Hamiltonian system has $\,n^2+n\,$ degrees of
freedom. There exist symplectic leafs of higher
codimensions (but always of an even dimension).

\medskip\noindent
{\bf Example: Case $n=1$.}
The two-system has order 5: two original variables $q$, $p$, and three ``moments
of the second order" $\alpha$, $\beta$, $\gamma$, cf.\ the system in the 
Introduction.
The explicit expression for the matrix $\Omega$ 
is given by
\begin{equation}
\label{omeg}
\Omega=J\oplus \omega,\quad
J=\left[\begin{array}{r r} 0&1\\-1&0\end{array}\right],\quad
\omega=\left[\begin{array}{r r r} 0& 2\al&4\be\\-2\al &0 &2\ga\\
-4\be&-2\ga &0 \end{array} \right].
\end{equation}
The (degenerate) Poisson bracket on the phase space is
$\;\{U,V\}=\sum_{i,j=1}^5 \Omega_{ij}\,\partial_i U\,
\partial_j V$.
The bracket has one Casimir function
\begin{equation}
\label{del}
\del=\al\be-{\ga}^2.
\end{equation}

\subsection{Integrable cases and special solutions }

Integrability of the two-system (\ref{eq8}) has a rare occurence.
We discuss two cases where it takes place.

\begin{Lemma} 
$\,$
 Two-system of a Hamiltonian system with Hamilton's
function being a polynomial of degree $2$ or less
splits into two independent subsystems and has the
general solution in a closed form.
\end{Lemma}

\noindent
{\it Proof}.
For any $i,j,k$, $\;\;H'''_{ijk}=0$,
and the matrix $H''$ of the second derivatives is constant.
Thus the system (\ref{eq8})  splits into the original system
(\ref{eq1}) and the system for $\Phi(t)$ with constant matrix $A$.
Both subsystems are linear, have constant coefficients, and therefore
integrable in elementary functions.
\WhiteBox

\begin{Lemma}
$\,$
Two-system of an integrable Hamiltonian system written in 
action-angle variables, with Hamilton's function depending
only on the action variables, is integrable.
\end{Lemma}

\noindent
{\it Proof}. We'll consider only case $n=1$, where the
nature of integrability is same as in the general case.
Let $(I,\theta)$ be the action-angle variables, and $\,H=H(I)\;$  Hamilton's
function. Denote
$M=\left(\begin{array}{l l} \al&\be\\
\be&\ga\end{array}\right)$.
The two-system reads
\begin{equation}
\label{eq11b}
 \dot I=0,\quad \dot \theta= H'(I)
 +H'''(I)\al,\quad
 \dot \al=0,\quad \dot\be=H''(I)\be,\quad
 \dot\ga =2 H''(I)\be
\end{equation}
The general solution of (\ref{eq11b}) is easily found to be 
$$
I=I_0,\quad \alpha=\alpha_0,\quad \theta=\theta_0 +t\omega_1,\quad
\beta=\beta_0 e^{\omega_2 t},
\quad
\gamma= 2\beta_0(e^{\omega_2t}-1)+\gamma_0.
$$
Here $\,I_0$, $\theta_0$, $\alpha_0$,
$\beta_0$, $\gamma_0\,$ are arbitrary constants and
$\;\omega_1=H'(I_0)+H'''(I_0)\alpha_0$, $\;\omega_2=H''(I_0)$.
\WhiteBox

\bigskip
Proposition 9 shows that the two-system in action-angle variables
is always integrable, while in general it is non-integrable except
in the situation of Proposition 8.
Thus the two-system is essentialy not invariant or covariant under canonical 
transformations of system (\ref{eq1}). 
Neveretheless,
there are two cases where the $\x$-component of a solution of the two-system
obeys equation (\ref{eq1}). 

\begin{Lemma} 
$\,$
{\rm(i)}$\,$ Consider special initial data for the two-system
$$
  \x(0)=\x_0,\quad \Phi(0)=0.
$$
The solution of {\rm (\ref{eq8})} with these initial data
is $(\x(t|\x_0),\;\Phi\equiv 0)$, where $\,\x(t|\x_0)\,$ is the solution
of {\rm (\ref{eq1})} with initial data $\x_0$.

\noindent
{\rm (ii)}$\,$ Let $x_0$ be such a stationary point of system (\ref{eq1})
where all third derivatives of the Hamiltonian  vanish:
$H'''_{ijk}(x_0)=0$.
Then the two-system {\rm (\ref{eq8})} has
a solution, for which $\,x(t)\equiv x_0\,$ and $\,\Phi(t)$
obeys the linear system with constant coefficients
\begin{equation}
\label{phi}
  \dot\Phi= [A(x_0),\Phi].
\end{equation}
\end{Lemma}

\noindent
{\it Proof} $\,$ easily follows from the structure of system (\ref{eq8}).
\WhiteBox


\end{document}